\documentclass[a4paper,draft,reqno,12pt]{amsart}
\usepackage[english]{babel}
\usepackage{amsmath}
\usepackage{amssymb}
\usepackage{amscd}
\usepackage{amsthm}
\usepackage{euscript}
\usepackage{tikz}
\usepackage{comment} 
\newtheorem{proposition}{Proposition}

\newtheorem{theorem}{Theorem}
\newtheorem*{theorem*}{Theorem}
\newtheorem{question}{Question}

\theoremstyle{definition}

\newtheorem{example}{Example}
\newtheorem{conjecture}{Conjecture}
\newtheorem{problem}{Problem}
\theoremstyle{remark}
\newtheorem {remark}{Remark}

\DeclareMathOperator{\Aut}{Aut}

\def\reg{{\mathrm{reg}}}

\def\BG{{\mathbb G}}
\def\BK{{\mathbb K}}

\def\BG{{\mathbb G}}

\def\BK{{\mathbb K}}

\def\BZ{{\mathbb Z}}

\def\BN{{\mathbb N}}
\def\BQ{{\mathbb Q}}

\def\CC{\mathcal{C}}

\def\CD{\mathcal{D}}

\def\CR{\mathcal{R}}

\def\SAut{\mathrm{SAut}}

\def\Cl{\mathrm{Cl}}

\def\div{\mathrm{div}}

\def\Xnod{{X}_{Y, B}}

\title[Flexibility of affine spherical varieties]{Flexibility of affine spherical varieties}

\thanks{The article was prepared within the framework of the project “International Academic Cooperation”
HSE University.
}

\author{Anton Shafarevich}
\email{shafarevich.a@gmail.com}
\address{
Lomonosov Moscow State University, Faculty of Mechanics and Mathematics, Department of Higher Algebra, Leninskie Gory 1, Moscow, 119991 Russia;
\linebreak
and
\linebreak
HSE University, Faculty of Computer Science, Pokrovsky Boulevard 11, Moscow, 109028, Russia}

\subjclass[2020]{Primary 14L30, 14J70; Secondary 13E10, 14L24.}
\keywords{Algebraic variety, algebraic group, action of additive group, spherical variety.}
\begin{document}
\maketitle

\begin{abstract}
We prove that the automorphism group $\Aut(X)$ of an affine spherical variety $X$ acts transitively on the set of smooth points $X^{reg}.$ If every invertible regular function on $X$ is constant, we prove that $X$ is flexible, i.e., the subgroup of $\Aut(X)$ generated by all $\BG_a$-subgroups acts transitively on $X^{reg}.$ 
 
 \end{abstract}

\section{Introduction}

Throughout this paper, $\BK$ denotes an algebraically closed field of characteristic zero. We denote by $\BG_a$ the group $(\BK,+).$ Every effective $\BG_a$-action on an algebraic variety $X$ corresponds to a subgroup $H$ in the automorphism group $\Aut(X)$. Such subgroups $H$ are called \emph{$\BG_a$-subgroups}. We denote by $\SAut(X)$ the subgroup of $\Aut(X)$ generated by all $\BG_a$-subgroups.

A smooth point $x$ of an irreducible algebraic variety $X$ is called \emph{flexible} if the tangent space $T_x X$ is spanned by the tangent vectors to the orbits of all $\mathbb{G}_a$-actions on $X$ that pass through $x$. The variety $X$ is called \emph{flexible} if every smooth point of $X$ is flexible. Affine flexible varieties have many beautiful properties. A key result in this context is the following theorem.

\begin{theorem}\cite[Theorem 0.1]{AFKKZ}\label{5autthm}
    Let $X$ be an affine irreducible variety of dimension at least two. Then the following conditions are equivalent.
    \begin{enumerate}
        \item $X$ is flexible;
        \item $\SAut(X)$ acts transitively on $X^{reg};$
        \item $\SAut(X)$ acts infinitely transitively on $X^{reg}$.
    \end{enumerate}
\end{theorem}

By infinite transitivity we mean that for every $m \in \mathbb{Z}_{>0}$ the action of $\SAut(X)$ on $X^{reg}$ is $m$-transitive, i.e., it is transitive on the set of ordered $m$-tuples of distinct points.

Note that the only one-dimensional affine variety admitting a non-trivial $\BG_a$-action is an affine line. Thus conditions 1 and 2 in Theorem \ref{5autthm} are still equivalent in dimension 1.

The simplest example of an affine flexible variety is the affine space $\mathbb{A}^n$. In \cite{AKZ}, the affine cones over flag varieties, as well as affine toric varieties and suspensions over affine spaces, were shown to be flexible. There are many results on  flexibility of affine cones over various projective varieties; see for example \cite{HT, KP, Pe}. Applications of the flexibility property can be found in \cite{Ar,Ka}. A broader review of flexibility can be found in \cite{Ar2}.

Every semisimple algebraic group is generated by its $\BG_a$-subgroups. Therefore, every affine homogeneous space of a semisimple algebraic group is flexible. The following problem was stated in \cite{AFKKZ2}.

\begin{problem}
    Characterize flexible varieties among affine varieties admitting an action of a semisimple algebraic group with a dense open orbit.
\end{problem}
So far, there are no known examples of non-flexible affine varieties admitting an action of a semisimple algebraic group with a dense open orbit.

A similar question arises when one considers reductive group actions. Throughout the paper, reductive groups are assumed to be connected. If $G$ is an algebraic group, then a \emph{$G$-variety} is an algebraic variety with a fixed regular action of $G$. We say that a $G$-variety $X$ is \emph{almost homogeneous} if there is an open dense $G$-orbit in $X$. It was proved in \cite{AFKKZ} that any affine, smooth, almost homogeneous $G$-variety for a semisimple group $G$ is flexible.

Any invertible regular function on an affine variety $X$ is invariant under all $\mathbb{G}_a$-actions. Consequently, if $X$ is flexible, then $\BK[X]^{\times} = \BK^{\times}$, where $\BK[X]$ denotes the algebra of regular functions on $X$ and $\BK[X]^{\times}$ is its group of invertible elements.  For a smooth affine $G$-variety $X$ of a reductive group $G$, the condition $\BK[X]^{\times} = \BK^{\times}$ implies flexibility; see \cite[Theorem~2]{GS}.

Among almost homogeneous varieties, there is an important class of spherical varieties. Let $G$ be a reductive group. An irreducible normal algebraic $G$-variety is called \emph{spherical} if a Borel subgroup of $G$ has an open orbit in $X$. The following conjecture was proposed in \cite[Section 6]{AFKKZ2}.

\begin{conjecture}\label{conj}
    Any affine spherical variety without non-constant invertible regular functions is flexible.  
\end{conjecture}

First, in \cite{AKZ}, Conjecture \ref{conj} was proved for toric varieties, that is, when the acting group $G$ is an algebraic torus. Subsequently, in \cite{Sh}, the conjecture was proved for horospherical varieties of semisimple groups. A horospherical variety is a $G$-variety such that the stabilizer of any point contains a maximal unipotent subgroup of $G$. Later, in \cite{GS}, the conjecture was proved for horospherical $G$-varieties with an arbitrary reductive group $G$. In this paper, we prove Conjecture \ref{conj} for all affine spherical varieties. 

\begin{theorem}\label{Thm2}
Let $G$ be a reductive group and $X$ an affine spherical $G$-variety. Suppose that $\BK[X]^{\times} = \BK^{\times}$. Then $X$ is flexible.  
\end{theorem}

For an arbitrary affine spherical $G$-variety, we prove the following result.

\begin{theorem}\label{Thm1}
Let $G$ be a reductive group and $X$ an affine spherical $G$-variety. Then the subgroup of $\Aut(X)$ generated by the image of $G$ in $\Aut(X)$ and all $\BG_a$-subgroups acts transitively on $X^{reg}.$

\end{theorem}

The class of affine spherical varieties contains many interesting examples. In addition to the affine toric and horospherical varieties mentioned above, this class includes affine cones over projective spherical varieties. Affine algebraic monoids with a reductive group of units are spherical varieties; see \cite{Vi, Ri}. The papers \cite{Pa1, Pa2} study spherical nilpotent orbits for the adjoint action of a semisimple group. Consequently, the normalization of the closure of any such orbit is also an affine spherical variety.

We now sketch the proofs of Theorems \ref{Thm2} and \ref{Thm1}. Let $G$ be a reductive group, $B\subseteq G$ a Borel subgroup, and $X$ a spherical $G$-variety. Denote by $\CD^B(X)$ the set of $B$-invariant prime divisors in $X$. For a $G$-orbit $Y$ in $X$, let $\CD^B_Y(X)$ be the set of $B$-invariant prime divisors in $X$ that contain $Y$. Then there is an open affine $B$-invariant subset

$$\Xnod = X \setminus \bigcup_{D\in \CD^B(X)\setminus\CD^B_Y(X)} D.$$

Let $P_Y$ be the subgroup of $G$ consisting of all elements that preserve $\Xnod$, and $P_u$ the unipotent radical of $P_Y$. By the Local Structure Theorem, there exist a Levi subgroup $L$ of $P_Y$ and a closed $L$-invariant subset $Z$ of $\Xnod$ such that $\Xnod \simeq P_Y \times^L Z \simeq P_u \times Z$; see Theorem \ref{LST}. The subvariety $Z$ is an affine spherical $L$-variety, and $Z \cap Y$ is the unique closed $L$-orbit in $Z$. If $Y$ is contained in $X^{reg}$, then $Z$ is smooth. In this case, the proof of Theorem 5.6 in \cite{AFKKZ} shows that the subgroup of $\operatorname{Aut}(Z)$ generated by $L$ and all $\BG_a$-subgroups acts transitively on $Z$. Therefore, the subgroup of $\operatorname{Aut}(\Xnod)$ generated by $P_Y$ and all $\BG_a$-subgroups acts transitively on $\Xnod$. For details, see Section~\ref{SecSph} and Proposition~\ref{LocalProp}.

If $X$ is affine, then $\Xnod$ is a principal open subset of $X$ by Proposition~\ref{Princ}. This implies that for any $\BG_a$-subgroup $H \subseteq \operatorname{Aut}(\Xnod)$ there exists a $\BG_a$-subgroup $H' \subseteq \operatorname{Aut}(X)$ such that for any point $x \in \Xnod$ the $H$-orbit of $x$ coincides with the $H'$-orbit of $x$; see Proposition~\ref{Ext}. Theorem \ref{Thm1} follows from Propositions \ref{Ext} and~\ref{LocalProp}.

Now suppose that $X$ is a normal affine variety with finitely generated divisor class group $\mathrm{Cl}(X)$, finitely generated Cox ring $\CR(X)$,  and $\BK[X]^{\times} = \BK^{\times}$. It was shown in \cite{GS} that if a subgroup of $\operatorname{Aut}(X)$ generated by all $\BG_a$-subgroups and a torus acts transitively on $X^{reg}$, then $X$ is flexible. We provide a simpler proof of this fact without restrictions on the class group and the Cox ring; see Proposition \ref{Prop1}. Finally, Proposition \ref{Prop1} implies Theorem \ref{Thm2}.   

\emph{Acknowledgment.} The author is grateful to Ivan Arzhantsev, Roman Avdeev, Vladimir Zhgoon and Dmitry Timashev for useful comments.

\section{Spherical varieties}\label{SecSph}

In this section, we recall some facts about spherical varieties. We refer to \cite{Kn, Perr} or \cite{Ti} for details. 

Let $G$ be a reductive group and $B$ a Borel subgroup in $G$. Denote by $\mathfrak{X}(B)$ the group of characters of $B$. It is a finitely generated free abelian group. Let $X$ be a spherical $G$-variety, that is, an irreducible normal $G$-variety which has an open $B$-orbit. For each character $\lambda \in \mathfrak{X}(B)$, there is a subspace $\BK(X)^{(B)}_{\lambda}$ in the field of rational functions $\BK(X)$ consisting of $B$‑semi‑invariants of weight $\lambda$. The set of characters for which the subspace $\BK(X)^{(B)}_{\lambda}$ is non-zero forms a lattice $M$, called the \emph{weight lattice} of $X$. We denote by $N$ the dual lattice $\mathrm{Hom}(M,\BZ)$, and by $M_{\BQ}$ and $N_{\BQ}$ the $\BQ$-vector spaces $M\otimes_{\BZ}\BQ$ and $N\otimes_{\BZ}\BQ$, respectively. There is a natural pairing 
$$\langle \cdot , \cdot \rangle:N\times M \to \mathbb{Z},\ \langle n, m\rangle = n(m). $$
This pairing can be extended to a pairing between $N_{\BQ}$ and $M_{\BQ}$.

Let us choose a point $x_o$ in the open $B$-orbit in $X$. Then for every $\lambda \in M$ there is a unique rational function $f_{\lambda} \in \BK(X)^{(B)}_{\lambda}$ such that $f_{\lambda}(x_o) = 1.$ Then we have $f_{\lambda + \mu} = f_{\lambda}\cdot f_{\mu}.$ Any $\BQ$-valued discrete valuation $v$ of the field $\BK(X)$ vanishing on $\BK^{\times}$ defines an element $\varphi(v) \in N_{\BQ}$ by the following rule:
$$\langle \varphi(v), \lambda\rangle = v(f_{\lambda}), $$
for all $\lambda \in M.$ The image of the set of  all $G$-invariant $\BQ$-valued discrete valuations under the map $v \to \varphi(v)$ is a finitely generated cosimplicial cone, called the \emph{valuation cone}. We denote it by~$\mathcal{V}$.

Consider the sets $\mathcal{D}^G(X)$ and $\mathcal{D}^B(X)$ of all $G$-invariant and $B$-invariant prime divisors on~$X$. Elements of the set $\mathcal{D} = \mathcal{D}^B(X)\setminus \mathcal{D}^G(X)$ are called \emph{colors}. For each divisor $D \in \CD^B(X)$, there is a discrete valuation $v_D$ of the function field $\BK(X)$ associated with $D$. This valuation, in turn, corresponds to a vector $\varphi(v_D)$ in the space $N_{\BQ}$. Thus, we obtain a map 
$$\rho: \mathcal{D}^B(X)\to N_{\BQ},\ D\to \rho(D) = \varphi(v_D).$$
The map $\rho$ is injective on $\CD^G(X).$

Now let $X$ be a simple spherical $G$-variety, that is, $X$ has only one closed $G$-orbit $O_c$. Note that all affine spherical $G$-varieties are simple. We define the set $\mathcal{D}_X = \{ D\in \mathcal{D} \mid O_c \subseteq D\}$ and the cone
$$\mathcal{C}_X = \mathrm{cone}(\rho(D) \mid D\in \CD^G(X)\sqcup \CD_X) \subseteq N_{\BQ}.$$
The cone $\CC_X$ is a pointed polyhedral cone. The pair $(\CC_X, \CD_X)$ is called the \emph{colored cone} of the spherical variety $X$. The map $D \to \rho(D)$ establishes a bijection  between the elements of $\mathcal{D}^G(X)$ and the primitive vectors on the edges of the cone $\CC_X$ that do not contain $\rho(D')$ for any  $D'\in \mathcal{D}_X$; see \cite[Lemma 3.4]{Kn}.

For every $G$-orbit $Y$ in $X$, consider the subset $X_Y = \bigsqcup Y'$, where $Y'$ runs over all $G$-orbits in $X$ such that $Y \subseteq \overline{Y'}$. The subset $X_Y$ is an open $G$-invariant subset of $X$ and $X_Y$ is a simple spherical $G$-variety with closed $G$-orbit $Y$. If $(\CC_Y, \CD_Y)$ is the colored cone of $X_Y$, then $\CC_Y$ is a face of $\CC_X$ and $\CD_Y = \CD_X\cap \rho^{-1}(\CC_Y).$ Denote by $\CD_Y^B(X)$ the set of $B$-invariant divisors in $X$ that contain the orbit $Y$. Then the set 
$$\Xnod = X \setminus \bigcup_{D\in \CD^B(X)\setminus\CD^B_Y(X)} D$$
is open, affine, $B$-invariant, and intersects $Y$; see \cite[Theorem 2.1]{Kn}.
\begin{proposition}\label{Princ}
    If $X$ is affine, then $\Xnod$ is a principal open subset.
\end{proposition}

\begin{proof}
    In fact, this statement follows from the proof of Theorem 3.1 in \cite{Kn}. However, we will provide a proof that uses only the combinatorics of spherical varieties. By Theorem 6.7 in \cite{Kn}, there is $\chi\in M$ with 
    $$\chi|_{\mathcal{V}} \leq 0,\ \chi|_{\CC_X} =  0,\ \chi|_{\rho(\mathcal{D}\setminus \CD_X)}  > 0.$$
    Consider the cone 
    $$\widehat{\CC_X} = \mathrm{cone}(\rho(D)\mid D\in \CD^B(X)).$$ 
    Then $\CC_X$ is a face of $\widehat{\CC_X}$ defined by the equation $\chi =0$. Therefore, $\CC_Y$ is a face of $\widehat{\CC_X}$ and there is a face $\CC_Y^{\perp}$ of the dual cone $\widehat{\CC_X}^{\vee}\subseteq M_{\BQ}:$
    $$\CC_Y^{\perp} = \{m\in \widehat{\CC_X}^{\vee}\mid m|_{\CC_Y} = 0\}.$$
    Consider $\lambda \in M$ that lies in the relative interior of the cone $\CC^{\perp}_Y.$ Let us take $B$-semi-invariant rational function $f_{\lambda}\in\BK(X)^{(B)}$ corresponding to $\lambda.$ Then 
    $$\div(f_{\lambda}) = \sum_{D\in \CD^B(X)} \langle  \rho(D), \lambda\rangle D. $$
    We have $\langle \rho(D),\lambda\rangle = 0$ for $D\in\CD^B(X)$ if and only if $\rho(D) \in \CC_Y.$ But $\rho(D)\in \CC_Y$ if and only if $D \in \CD_Y^B(X).$ For $D \in \CD^B(X)\setminus \CD_Y^B(X)$ we have $\langle \rho(D), \lambda\rangle>0.$ Therefore, $f_{\lambda}$ is a regular function on $X$ and 
    $$\{x\in X \mid f_{\lambda}(x) = 0\} = \bigcup_{D\in \CD^B(X)\setminus \CD^B_Y(X)} D.$$
    So $\Xnod = X_{f_{\lambda}} = \{x\in X \mid f_{\lambda}(x)\neq 0\}.$
\end{proof}

The structure of $\Xnod$ is described by the Local Structure Theorem. Let $P_Y$ be the parabolic subgroup of $G$ consisting of all elements that preserve $\Xnod$. Denote by $P_u$ the unipotent radical of $P_Y$. Choose a Levi subgroup $L$ in $P_Y$.

\begin{theorem}(\cite[Theorem 15.17]{Ti}; see also Section 1 in \cite{Br}.)\label{LST} There is an $L$-invariant closed subvariety $Z \subseteq \Xnod$ with the following properties:
\begin{enumerate}
    \item $\Xnod \simeq P_Y \times^L Z \simeq P_u\times Z;$
    \item $Z$ is a spherical $L$-variety and $Z\cap Y$ is a unique closed $L$-orbit in $Z$.
\end{enumerate}

\end{theorem}

Here by $P_Y\times^LZ$ we mean the total space of the homogeneous fiber bundle over $P_Y/L$ with fiber $Z$. By definition, it is the geometric quotient of $P_Y \times Z $ by the following action of $L$:
$$g\circ(p,z) = (pg^{-1},gz), \ \text{for}\ g\in L,\ (p,z) \in P_Y\times Z.$$
For details, see \cite[Section 4.8]{PV}. 

If the orbit $Y$ is smooth, then the variety $Z$ is also smooth. Let $H\subseteq L$ be the stabilizer of a point in the closed $L$-orbit $Z\cap Y$. Then Luna's Étale Slice Theorem implies that there is a finite-dimensional $H$-module $W$ such that the variety $Z$ is $L$-equivariantly isomorphic to the total space $L\times^HW$ of the homogeneous vector bundle over $L/H$ with the fiber $W$; see Theorem 6.7 in \cite{PV}. 

\begin{proposition}
    The subgroup in $\operatorname{Aut}(Z)$ generated by $L$ and all $\BG_a$-subgroups acts transitively on $Z$. 
\end{proposition}
\begin{proof}
The argument follows the proof of Theorem 5.6 in \cite{AFKKZ}. Since $L/H$ is an affine variety and $Z = L \times^H W$ is a vector bundle over $L/H$, for any $w\in W$, there is a section
$$s\in H^0(L/H, L\times^H W)$$ such that $s(eH) = w;$ see \cite[Example 5.16.2]{Ha}. Then one can define the following $\BG_a$-action on $Z$:
$$t\circ z = z + ts(\pi(z)),\ t\in \BG_a, z\in Z.$$
Here $\pi \to L/H$ denotes the bundle projection. Note that $z$ and $s(\pi(z))$ belong to the same fiber of $Z$, so their linear combination is well-defined. It follows that the subgroup of $\Aut(Z)$ generated by all $\BG_a$-actions acts transitively on each fiber of $Z$.

Since $L$ acts transitively on $L/H$, it follows that the group generated by $L$ and all $\BG_a$-subgroups acts transitively on $Z$.

\end{proof}

The group $V$ is isomorphic to  $\BG_a^m$, where $m = \dim V$ and the action of $V$ on $Z$ can be extended to an action of $V$ on $\Xnod \simeq P_u \times Z.$ So we obtain the following proposition.

\begin{proposition}\label{LocalProp}
    Let $X$ be an affine spherical $G$-variety and $Y\subseteq X^{reg}$ a $G$-orbit. Then the subgroup of $\Aut(\Xnod) $ generated by the image of $P_Y$ in $\Aut(\Xnod)$ and all $\BG_a$-subgroups acts transitively on~$\Xnod$.
\end{proposition}

\begin{remark}
    In fact, for a smooth orbit $Y$, the subset $\Xnod$ is isomorphic to a product of a torus and an affine space; see \cite[Corollary 3]{Br}. Proposition \ref{LocalProp} can be deduced from this.
\end{remark}

\section{$\BG_a$-actions on affine varieties}

Here we recall some facts about $\BG_a$-actions on affine varieties and prove Theorem~\ref{Thm1}. Let $X$ be an irreducible affine variety. A derivation of the algebra $\BK[X]$ is a $\BK$-linear map $\delta:\BK[X]\to\BK[X]$ satisfying the Leibniz rule. A derivation $\delta$ of the algebra $\BK[X]$ is called \emph{locally nilpotent} if for any $a\in \BK[X]$ there exists $n = n(a)\in \BN$ such that $\delta^n(a) = 0$. Given a locally nilpotent derivation of the algebra $\BK[X]$, one can construct a $\BG_a$-action on $\BK[X]$ by the following rule:

$$s\circ f = \exp(s\delta)(f) = \sum_{i\geq 0} \frac{s^i\delta^i(f)}{i!},\ s\in \BG_a,\ f\in \BK[X].$$
This $\BG_a$-action on $\BK[X]$ defines a $\BG_a$-action on $X$. Conversely, any $\BG_a$-action on $X$ is obtained via this construction; see Section 1.5 in \cite{Fr} for details.

\begin{proposition}\label{Ext}
    Let $X$ be an affine irreducible variety and $U\subseteq X$ a principal open subset. Then for any $\BG_a$-subgroup $H$ in $\Aut(U)$ there exists a $\BG_a$-subgroup $H'$ in $\Aut(X)$ such that $Hx = H'x$ for every $x\in U.$
\end{proposition}
\begin{proof}
    There is $h \in \BK[X]$ such that $U = \{x\in X \mid h(x) \neq 0\}$ and $\BK[U] = \BK[X]_h.$ The $\BG_a$-subgroup $H$ defines a $\BG_a$-action on $U.$ Let $\delta$ be the corresponding locally nilpotent derivation on $\BK[U]$. Since $h$ is invertible in $\BK[U]$, we have $\delta(h^k) = 0$ for all $k \in \BZ_{>0}$;  see \cite[Corollary 1.23]{Fr}. Moreover, the linear operator $h^k\delta$ is a locally nilpotent derivation on $\BK[U]$ by \cite[Principle 7]{Fr}.

    Let $f_1,\ldots, f_s$ be generators of $\BK[X] \subseteq \BK[U].$ Then $\delta(f_i) = \frac{g_i}{h^{m_i}}$ for some $g_i \in \BK[X]$ and $m_i \in \BN.$ Let $m$ be the maximum among $m_1,\ldots, m_s.$ Then $h^m\delta(f_i) \in \BK[X]$ for all $i=1\ldots, s.$ So $\BK[X]$ is an invariant subalgebra for $h^m\delta$ and $h^m\delta$ is a locally nilpotent derivation of $\BK[X].$ Let $H'$ be the corresponding $\BG_a$-subgroup in $\Aut(X).$

    For $x\in U, f\in \BK[X]$, and $s\in H$, the value of $f$ at $s\circ x$ is equal to
    $$f(s\circ x) = ((-s)\circ x)(f) = \sum_i \frac{(-s\delta)^i(f)(x)}{i!}.$$
    Now, if we act on $x$ by the element $\frac{s}{h^m(x)} \in H'$, the result is
    $$f\left(\frac{s}{h^m(x)}\circ x\right) = \left(-\frac{s}{h^m(x)}\circ f\right)(x) = \sum_i \frac{(-s\delta)^i(f)(x)}{i!}.$$
    So the image of $x$ with respect to $s\in H$ is equal to the image of $x$ with respect to $\frac{s}{h^m(x)}\in H'.$
    
\end{proof}

Now we are ready to prove Theorem \ref{Thm1}. 

\begin{proof}[Proof of Theorem \ref{Thm1}]
    Denote by $\mathrm{GAut}(X)$ the subgroup of $\Aut(X)$ generated by $G$ and all $\BG_a$-subgroups of $\Aut(X)$. Let $Y\subseteq X$ be a $G$-orbit contained in $X^{reg}$. By Proposition \ref{LocalProp}, the subgroup of $\Aut(\Xnod)$ generated by $P_Y$ and all $\BG_a$-subgroups of $\Aut(\Xnod)$ acts transitively on $\Xnod.$ By Proposition \ref{Princ}, $\Xnod$ is a principal open subset of $X$. Applying Proposition \ref{Ext}, we obtain that $\Xnod$ is contained in a single $\mathrm{GAut}(X)$-orbit. Since $\Xnod$ intersects both the $G$-orbit $Y$ and the open $G$-orbit $O$, it follows that $Y$ lies in the same $\mathrm{GAut}(X)$-orbit as $O$.
\end{proof}

\section{Flexibility}

In this section we discuss the flexibility property. Let us recall some facts about flexible points and the subgroup $\SAut(X)$.

\begin{proposition}[Proposition 1.3, Corollary 1.11 and Theorem 1.13 in \cite{AFKKZ}]\label{Flexprop} Let $X$ be an irreducible algebraic variety. Then
\begin{enumerate}
    \item the orbits of $\SAut(X)$ in $X$ are locally closed;
    \item there is a finite set of rational $\SAut(X)$-invariants which separate orbits in some nonempty open subset of $X$;
    \item a point $x\in X^{reg}$ is flexible if and only if the $\SAut(X)$-orbit of $x$ is open in $X$.
\end{enumerate}
    
\end{proposition}

Suppose $X$ is a normal affine variety such that both the divisor class group $\Cl(X)$ and the Cox ring $\CR(X)$ are finitely generated, and $\BK[X]^{\times} = \BK^{\times}$. By \cite[Lemma 5]{GS}, if the subgroup of $\Aut(X)$ generated by all $\BG_a$-subgroups and a torus acts transitively on $X^{reg}$, then $X$ is flexible. We give a simpler proof that does not require the assumptions on the class group or the Cox ring.

\begin{proposition}\label{Prop1}
    Let $G$ be a connected linear algebraic group acting effectively on an  algebraic variety $X$ (not necessarlu affine). Assume that the group $\langle G, \SAut(X) \rangle$ acts transitively on $X^{reg}$ and $\BK[X^{reg}]^{\times} = \BK^{\times}$. Then $\SAut(X)$ acts transitively on $X^{reg}.$
\end{proposition}

\begin{proof}
    This proof was communicated to the author by Dmitriy Timashev.

    We have $G = G_{red} \ltimes R_u$ where $G_{red}$ is a reductive subgroup and $R_u$ is the unipotent radical. Moreover, $G_{red}$ is an almost direct product of the derived subgroup $G_s = [G_{red}, G_{red}]$ and the radical $T$, which is a torus. The subgroup $G_s$ is a semisimple group. Thus $G_s$ is generated by its $\BG_a$-subgroups. The group $R_u$ is unipotent and is also generated by $\BG_a$-subgroups. So
    $$\langle G, \SAut(X)\rangle = \langle T, \SAut(X) \rangle.$$

    Suppose that the action of $\SAut(X)$ on $X^{reg}$ is not transitive. The orbits of $\SAut(X)$ are locally closed. Let $Y$ be an $\SAut(X)$-orbit in $X^{reg}$. Consider a subgroup 
    $$S = \{t\in T\mid tY = Y\}.$$
    It is easy to check that the map 
    $$\varphi: T\times^SY \longrightarrow X^{reg},\ [(t,y)] \longrightarrow t\circ y$$
    is well-defined and an isomorphism. Here, by $[(t,y)]$ we mean the $S$-orbit of the element $(t,y)\in T\times Y$.

    If $S \neq T$, then there is a non-constant invertible regular function $f$ on $T/S$. If we denote by $\pi: T\times^SY \to T/S$ the fiber projection, the function $\pi^*(f)$ is a non-constant invertible function on $T\times^S Y \simeq X^{reg}.$ So $T = S$ and $Y = X^{reg}.$

\end{proof}

\begin{proof}[Proof of Theorem \ref{Thm2}]
    Since $X$ is normal $\BK[X] = \BK[X^{reg}].$ By Theorem \ref{Thm1} the group $\langle G, \SAut(X)\rangle$ acts on $X^{reg}$ transitively. By Proposition \ref{Prop1} we obtain that if $\BK[X]^{\times} = \BK^{\times}$ then $X$ is flexible. 
\end{proof}

\begin{example}

We thank the referee of an earlier version of this paper for pointing out that Theorem \ref{Thm2} implies that every normal affine cone $X$ over a projective spherical variety $Y$ is flexible. Indeed, since $X$ is an affine cone, its coordinate ring admits a positive grading

$$
\BK[X]= \bigoplus A_n,$$
where $A_n = H^0(Y, L^n)$ for some ample line bundle $L$ on $Y$. We have $A_0 = \BK[Y] = \BK.$
Let $f\in \BK[X]^\times$, and write
$$
f=f_0+\cdots+f_m,\qquad f_i\in A_i,\quad f_m\neq 0.
$$
Similarly, write $f^{-1}=g_0+\cdots+g_n$, with $g_n\neq 0$. Since $X$ is normal and irreducible, $\BK[X]$ is a domain. The homogeneous component of degree $m+n$ in the equality $ff^{-1}=1$ is
$$
f_mg_n=0,
$$
which is impossible. Hence $m=n=0$, and therefore $f\in A_0=\BK$. Thus
$$
\BK[X]^\times=\BK^\times.
$$

Since $Y$ is spherical under the action of a reductive group $G$, the affine cone $X$ is spherical under the induced action of $G\times\BK^\times$, with $\BK^\times$ acting by scalar multiplication.
\end{example}

Note that, in general, the existence of an open SAut(X)-orbit on X does not imply flexibility. Examples can be found in \cite{Ko} and \cite{Ga}.

In the definition of a spherical variety, the normality condition is sometimes omitted. Flexibility of non-normal toric and horospherical varieties has been studied in \cite{BG} and \cite{GK}. 

Let $F$ be a subgroup of $\Aut(X)$ and $H$ a $\BG_a$-subgroup. Then $H$ is called an $F$-root subgroup if $F$ normalizes $H$. Let $T$ be an algebraic torus. The proof of flexibility of affine toric $T$-varieties in \cite{AKZ} used the description of $T$-root subgroups obtained in \cite{De} and \cite{Li}. The analog of $T$-root subgroups in the spherical case is the notion of $B$-root subgroups. Recently, significant progress has been made in the study of $B$-root subgroups on spherical varieties; see \cite{AA, AZ1, AZ2}.

\begin{question}
    Let $G$ be a reductive group and $X$ an affine spherical $G$-variety. Is it true that the subgroup in $\Aut(X)$ generated by the image of $G$ in $\Aut(X)$ and all $B$-root subgroups acts transitively on $X^{reg}$?
\end{question}

After the first version of this paper was posted on arXiv, the question was answered affirmatively in~\cite{AZ3}.

\begin{theorem}\cite[Corollary 1.2]{AZ3}\label{AZ3}
   Let $G$ be a reductive group and $X$ be a quasiaffine spherical $G$-variety (not necessarily normal). Let $\widetilde{G}$ be the subgroup in $\Aut(X)$ generated by the images of $G$ and all $B$-root subgroups on $X$. Then $X^{reg}$ is a single $\widetilde{G}$-orbit.

\end{theorem} 

Moreover, we can give a criterion for flexibility of non-normal quasiaffine spherical varieties.

\begin{theorem}
    Let $G$ be a reductive group and $X$ a quasiaffine spherical $G$-variety. Then $\SAut(X)$ acts transitively on $X^{reg}$ if and only if $\BK[X^{\reg}]^{\times} = \BK^{\times}.$
\end{theorem}
\begin{proof}
    If $\BK[X^{reg}]^\times = \BK^\times$, then by Theorem \ref{AZ3} and Proposition \ref{Prop1} the group $\SAut(X)$ acts transitively on $X^{reg}.$  If $\BK[X^{reg}]^\times \neq \BK^\times$, then any non-constant regular function $f \in \BK[X^{reg}]^\times $ is invariant under $\SAut(X)$. Thus  $\SAut(X)$ does not act transitively on $X^{reg}$.

\end{proof}

For affine toric varieties, all $\Aut(X)$-orbits were described in \cite{AB}. For horospherical varieties, a partial result was obtained in \cite{BGS}.

\begin{problem}
    Describe all $\Aut(X)$-orbits on affine spherical varieties.
\end{problem}


\begin{thebibliography}{99}

\bibitem{Ar2} I.~Arzhantsev. \textit{Automorphisms of algebraic varieties and infinite transitivity}. St. Petersburg Math. J. \textbf{34} (2023), no.~2, 143--178.

\bibitem{Ar} I.~Arzhantsev. \textit{On images of affine spaces}. Indag. Math. \textbf{34} (2023), no.~4, 812--819.

\bibitem{AA} I.~Arzhantsev and R.~Avdeev. \textit{Root subgroups on affine spherical varieties}. Selecta Math. (N. S.) \textbf{28} (2022), no.~3, Art.~60, 37~pp.

\bibitem{AB} I.~Arzhantsev and I.~Bazhov. \textit{On orbits of the automorphism group on an affine toric variety}. Cent. Eur. J. Math. \textbf{11} (2013), no.~10, 1713--1724. 

\bibitem{AFKKZ} I.~Arzhantsev, H.~Flenner, S.~Kaliman, F.~Kutzschebauch, and M.~Zaidenberg. \textit{Flexible varieties and automorphism groups}. Duke Math. J. \textbf{162} (2013), no.~4, 767--823.

\bibitem{AFKKZ2} I.~Arzhantsev, H.~Flenner, S.~Kaliman, F.~Kutzschebauch, and M.~Zaidenberg. \textit{Infinite transitivity on affine varieties}. In: Birational geometry, rational curves, and arithmetic, Simons Symp., Springer, Cham, 2013, pp.~1--13.

\bibitem{AKZ} I.~Arzhantsev, K.~Kuyumzhiyan, and M.~Zaidenberg. \textit{Flag varieties, toric varieties, and suspensions: three instances of infinite transitivity}. Sb. Math. \textbf{203} (2012), no.~7, 923--949.

\bibitem{AZS} I.~Arzhantsev, Y.~Zaitseva, and K.~Shakhmatov. \textit{Homogeneous Algebraic Varieties
and Transitivity Degree}.  Proc. Steklov Inst. Math. \textbf{318} (2022), 13--25.

\bibitem{AZ1} R.~Avdeev and V.~Zhgoon. \textit{On the existence of B-root subgroups on affine spherical varieties}. Dokl. Math. \textbf{105} (2022), no.~2, 51--55.

\bibitem{AZ2} R.~Avdeev and V.~Zhgoon. \textit{Root subgroups on horospherical varieties}. arXiv:2312.03377 (2023).


\bibitem{AZ3} R.~Avdeev and V.~Zhgoon. \textit{Connecting orbits in quasiaffine spherical varieties via $B$-root subgroups}. arXiv:2512.09906v2 (2026).


\bibitem{BGS} V.~Borovik, S.~Gaifullin, and A. Shafarevich. \textit{On orbits of automorphism groups on horospherical varieties.} Math. Nachr. \textbf{297} (2024), no. 9, 3174--3183.

\bibitem{BHM} A.~Bialynicki-Birula, G.~Hochschild, and G.~Mostow.
\textit{Extensions of Representations of Algebraic Linear Groups}.
Amer.~J.~Math. \textbf{85} (1963), no.~1, 131--144.

\bibitem{BG} I.~Boldyrev and S.~Gaifullin. \textit{Automorphisms of nonnormal toric varieties}. Math. Notes \textbf{110} (2021), no.~5, 872--886.

\bibitem{Br} M.~Brion. \textit{Groupe de Picard et nombres caractéristiques des variétés sphériques}. Duke Math. J. \textbf{58} (1989), no.~2, 397--424.


\bibitem{De} M.~Demazure. \textit{Sous-groupes algébriques de rang maximum du groupe de Cremona}. Ann. Sci. École Norm. Sup. \textbf{4} (1970), no.~3, 507--588.

\bibitem{FKZ} H.~Flenner, S.~Kaliman and M.~Zaidenberg.
\textit{A Gromov--Winkelmann type theorem for flexible varieties.}
J. Eur. Math. Soc. \textbf{18} (2016), no.~11, 2483--2510.

\bibitem{Ga} S.~Gaifullin. \textit{Generically flexible affine varieties with invariant divisors}. arXiv:2507.14745 (2025). 

\bibitem{GK} S.~Gaifullin and V.~Kikteva. \textit{Flexibility criterion for affine horospherical varieties}. Results Math. \textbf{81} (2026), no.~5, Art.~146.

\bibitem{GS} S.~Gaifullin and A.~Shafarevich. \textit{Flexibility of normal affine horospherical varieties}. Proc. Amer. Math. Soc. \textbf{147} (2019), no.~8, 3317--3330.

\bibitem{Fr} G.~Freudenburg. \textit{Algebraic theory of locally nilpotent derivations}. Second edition, Encyclopaedia Math. Sci., vol.~136, Springer, Berlin, 2017.

\bibitem{HT} N.~Hang and H.~Truong. \textit{The affine cones over Fano-Mukai fourfolds of genus 7 are flexible}. Int. Math. Res. Not. IMRN (2024), no.~10, 8417--8426.

\bibitem{Ha} R.~Hartshorne. \textit{Algebraic Geometry}. Graduate Texts in Mathematics, vol.~52, Springer, New York, 1977.


\bibitem{Ka} S.~Kaliman. \textit{Embedding theorems for flexible varieties}. Michigan Math. J. \textbf{75} (2025), no.~3, 601--614.



\bibitem{Kn} F.~Knop. \textit{The Luna-Vust theory of spherical embeddings}. In: Proceedings of the Hyderabad Conference on Algebraic Groups (Hyderabad, 1989), 225--249, Manoj Prakashan, Madras, 1991.

\bibitem{KKLV} F.~Knop, H.~Kraft, D.~Luna, and Th.~Vust. \textit{Local properties of algebraic group actions}. In: Algebraische Transformations gruppen und Invarianten theorie, DMV Seminar, vol.~13, pp.~63--75, Birkhäuser, Basel, 1989.

\bibitem{Ko} S.~Kovalenko. \textit{Transitivity of Automorphism Groups of Gizatullin Surfaces}. Int. Math. Res. Not. IMRN (2015), no.~21, 1433--1484.

\bibitem{KP} J.~Kim and J.~Park. \textit{Generic flexibility of affine cones over del Pezzo surfaces of degree 2}. Internat. J. Math. \textbf{32} (2021), no.~14, Art.~2150104, 18~pp.


\bibitem{Li} A.~Liendo. \textit{Affine T-varieties of complexity one and locally nilpotent derivations}. Transform. Groups \textbf{15} (2010), no.~2, 389--425.

\bibitem{Pa1} D. Panyushev. \textit{Complexity and nilpotent orbits}. Manuscripta Math. \textbf{83} (1994),
223--237.

\bibitem{Pa2} D. Panyushev. \textit{On spherical nilpotent orbits and beyond}. Annales l’Institut Fourier \textbf{49} (1999), 1453--1476.

\bibitem{Pe} A.~Perepechko. \textit{Flexibility of affine cones over del Pezzo surfaces of degree 4 and 5}. Funktsional. Anal. i Prilozhen. Funct. Anal. Appl. \textbf{47} (2013), no.~4, 284--289.

\bibitem{Perr} N.~Perrin. \textit{On the geometry of spherical varieties}. Transform. Groups \textbf{19} (2014), no.~1, 171--223.

\bibitem{PV} V.~Popov and E.~Vinberg. \textit{Invariant theory}. In: Algebraic Geometry. IV: Linear Algebraic Groups, Invariant Theory, Encycl. Math. Sci., vol.~55, pp.~123--278, Springer, Berlin, 1994.

\bibitem{Ri} A.~Rittatore. \textit{Algebraic monoids with affine unit group are affine}. Transform. Groups \textbf{12} (2007),
no.~3, 601--605.

\bibitem{Sh} A.~Shafarevich. \textit{Flexibility of S-varieties of semisimple groups}. Sb. Math. \textbf{208} (2017), no.~2, 285--310.

\bibitem{Ti} D.~Timashev. \textit{Homogeneous spaces and equivariant embeddings}. Encycl. Math. Sci., vol.~138, Springer, Berlin, 2011.

\bibitem{Vi} E.~Vinberg. \textit{On reductive algebraic semigroups}. In: Lie Groups and Lie Algebras, E.B. Dynkin
Seminar, S. Gindikin and E. Vinberg, Editors, Amer. Math. Soc. Transl. vol.~169, pp.~145--182, Amer. Math.
Soc., 1995.
\end{thebibliography}
\end{document}